\documentclass{amsart}

\usepackage{amscd,amssymb,amsmath}
\usepackage{graphics}
\usepackage[all]{xy}

\newtheorem{theorem}{Theorem}[section]
\newtheorem*{maintheorem}{Theorem}
\newtheorem{lemma}[theorem]{Lemma}
\newtheorem{proposition}[theorem]{Proposition}
\newtheorem{corollary}[theorem]{Corollary}

\theoremstyle{definition}

\newtheorem{remark}[theorem]{Remark}
\newtheorem*{acknowledgement}{Acknowledgement}

\theoremstyle{remark}

\renewcommand{\labelenumi}{(\roman{enumi})}

\DeclareFontFamily{U}{wncy}{}
\DeclareFontShape{U}{wncy}{m}{n}{<->wncyr10}{}
\DeclareSymbolFont{mcy}{U}{wncy}{m}{n}
\DeclareMathSymbol{\Sh}{\mathord}{mcy}{"58}

\newcommand\mynote[1]{\mbox{}\marginpar{\ \\ \small \tt #1}}
\newcommand\bel[1]{{\mynote{#1}}\begin{equation}\label{#1}}
\newcommand\mylabel[1]{\label{#1}}

\newcommand{\ZZ}{\mathbb{Z}}

\newcommand{\FF}{\mathbb{F}}

\newcommand{\PP}{\mathbb{P}}

\newcommand  {\shH}     {\mathcal{H}}

\newcommand  {\shK}     {\mathcal{K}}

\newcommand  {\shL}     {\mathcal{L}}

\newcommand  {\GL}      {\operatorname{GL}}

\newcommand  {\length}  {\operatorname{length}}

\newcommand  {\lra}     {\longrightarrow}

\newcommand  {\Mat}     {\operatorname{Mat}}
\newcommand  {\maxid}   {\mathfrak{m}}

\newcommand  {\NS}      {\operatorname{NS}}

\renewcommand{\O}       {\mathcal{O}}

\newcommand  {\ord}     {\operatorname{ord}}

\newcommand  {\Pic}     {\operatorname{Pic}}

\newcommand  {\pr}      {\operatorname{pr}}
\newcommand  {\Proj}    {\operatorname{Proj}}

\newcommand  {\quadand} {\quad\text{and}\quad}

\newcommand  {\ra}      {\rightarrow}

\newcommand  {\red}     {{\operatorname{red}}}

\newcommand  {\Sing}    {\operatorname{Sing}}

\newcommand  {\Spec}    {\operatorname{Spec}}

\def\mydate{\number\day\space\ifcase\month \or January\or February\or March\or 
April\or May\or June\or July\or
August\or September\or October\or November\or December\fi \space\number\year}

\DeclareFontFamily{U}{wncy}{}
\DeclareFontShape{U}{wncy}{m}{n}{<->wncyr10}{}
\DeclareSymbolFont{mcy}{U}{wncy}{m}{n}
\DeclareMathSymbol{\Sh}{\mathord}{mcy}{"58}

\begin{document}

\title[Surfaces of General Type]
      {Wildly Ramified Actions and Surfaces of General Type Arising from Artin--Schreier Curves}

\author[Hiroyuki Ito]{Hiroyuki Ito}
\address{Department of Mathematics,
Faculty of Science and Technology,
Tokyo University of Science,
2641 Yamazaki, Noda, Chiba, 278-8510, Japan}
\curraddr{}
\email{ito\_hiroyuki@ma.noda.tus.ac.jp}

\author[Stefan Schr\"oer]{Stefan Schr\"oer}
\address{Mathematisches Institut, Heinrich-Heine-Universit\"at,
40204 D\"usseldorf, Germany}
\curraddr{}
\email{schroeer@math.uni-duesseldorf.de}


\dedicatory{Dedicated to Gerard van der Geer\\Revised version, 20 February 2012}

\begin{abstract}
We analyse the  diagonal quotient for the product of certain Artin--Schreier curves.
The smooth models are almost always surfaces of general type, with Chern slopes
tending asymptotically to 1.
The calculation of   numerical invariants
relies on a close examination of the relevant wild quotient singularity in characteristic $p$.
It turns out that the canonical model has $q-1$ rational double points
of type $A_{q-1}$, and embeds as a divisor of degree $q$ in $\PP^3$,
which is in some sense reminiscent of the classical Kummer quartic.
\end{abstract}

\maketitle
\tableofcontents
\renewcommand{\labelenumi}{(\roman{enumi})}

\section*{Introduction}

It is a classical fact in complex geometry   that the singular Kummer surface $A/\left\{\pm 1\right\}$
attached to an abelian surface $A$ has sixteen rational double points, and an irreducible principal polarization
embeds it as a quartic surface in $\PP^3$, compare Hudson's classical monography \cite{Hudson 1903}, 
or for a modern account \cite{Gonzalez-Dorrego 1994}.
Our starting point was an analogous computation
in characteristic $p=3$ for the diagonal action  
of the additive group $G=\FF_3$ on the selfproduct $A=E\times E$, where $E:\, y^2=x^3+x$
is the supersingular elliptic curve, viewed as an Artin--Schreier covering.
It turns out that this is a special case of a  rather  general construction, which works
for all primes $p$, in fact for all prime powers $q=p^s$. It starts
with certain \emph{Artin--Schreier curves}  and   leads, with a few exceptions for small prime powers, to    
\emph{surfaces of general type}.

The goal of this paper is to describe the geometry of these surfaces, and we obtain a fairly complete description. 
The construction goes as follows: Fix an algebraically closed ground field $k$ of characteristic $p>0$ and consider  
Artin--Schreier curves of the form
$$
C:\quad f(y) = x^q-x,
$$
where $f$ is  a monic polynomial of degree $\deg(f)=q-1$. These curves
carry a translation action of the additive group $G=\FF_q\simeq (\ZZ/p\ZZ)^{\oplus s}$, and we consider
the diagonal action on the product $C\times C'$ of two such Artin--Schreier curves.
The quotient $(C\times C')/G$ is a normal surface containing a unique singularity.

Such singularities are examples of \emph{wild  quotient singularities}, i.e., the characteristic of the ground field
divides the order of the   group $G$. Few examples of wild quotient singularities occur in the literature, and
little is known in general.   Artin \cite{Artin 1975} gave a complete classification for wild  $\ZZ/2\ZZ$-quotient singularities
in dimension two,
and  general $\ZZ/p\ZZ$-quotient singularities were studied further by Peskin \cite{Peskin 1983}. Peculiar properties
of wild $S_n$-quotient singularities in relation to punctual Hilbert schemes appear in \cite{Schroeer 2009}.
In light of the scarcity of examples, it is useful to have more classes of wild quotient singularities in which
computations are feasible.

Lorenzini initiated a general investigation of wild quotient singularities on surfaces \cite{Lorenzini 2006}, which play
an important role in understanding the reduction behaviour of curves over discrete valuation fields.
He   compiled a list of open questions  \cite{Lorenzini 2011a}.
In a recent paper, Lorenzini  studied wild $\ZZ/p\ZZ$-quotient singularities resulting from diagonal 
actions on products of $C\times C'$, where one or both factors are ordinary curves \cite{Lorenzini 2011b}. In some sense, we treat the opposite
situation, as Artin-Schreier curves have vanishing $p$-rank, and our main concern is the interplay between
the local structure coming from the the wild quotient singuarity and the global geometry of the algebraic surface.
We have choosen   Artin--Schreier coverings that are concrete enough so that explicit computations are   possible.
Note that our set-up includes  wild quotient singularities with respect to elementary abelian groups, and not ony cyclic groups.

Consider the minimal resolution of singularities
$$
X\lra(C\times C')/G
$$
of our normal surface. According to Lorenzini's  general observation,
the exceptional divisor of a wild quotient singularity in dimension two consists of projective lines, and its
dual graph is a tree  \cite{Lorenzini 2006}. Using an explicit formal equation for the singularitiy, we show
that   the dual graph is even star-shaped, with $q+1$ terminal chains   attached to the central node,
each of length $q$, as depicted in Figure 1. 
The two basic numerical invariants of surface singularities are the genus $p_f=h^1(\O_Z)$ of the fundamental cycle $Z\subset X$
and the \emph{geometric genus} $p_g=h^1(\O_{nZ})$, $n\gg 0$. The latter is usually very difficult to compute.
We obtain 
$$
p_f=(q-1)(q-2)/2\quadand p_g=q(q-1)(q-2)/6,
$$
the latter under the assumption that our prime power $q$ is prime.
This relies on a computation of the global $l$-adic Euler characteristic of the surface, combined with a determination 
of the $G$-invariant part in $H^0(C\times C',\Omega^1_{C\times C'})$, which in 
turn depends on a problem in modular representation theory related to   tensor products of Jordan matrices. 
In contrast, the singularities occuring in \cite{Lorenzini 2011b}, where at least one factor in $C\times C'$ is an ordinary curve,
are all rational.

Having a good hold on the structure of the resolution of singularities, we   determine the global invariants
of the smooth surface $X$. It turns out that $H^1(X,\O_X)=0$, and that $\Pic(X)$ is a free abelian group of finite rank.
Moreover, the algebraic fundamental group $\pi_1(X)$ vanishes. Passing to the minimal model
$$
X\lra S,
$$
we show that the minimal surface $S$ is of general type for $q\geq 5$, a K3 surface for $q=4$, and
a weak del Pezzo surface for $q=2,3$. Their Chern invariants are given by the formula
$$
c_1^2=q^3-8q^2+16q\quad\quadand c_2=q^3-4q^2+6q \qquad\quad (q\geq 4).
$$
The resulting \emph{Chern slopes}  asymptotically  tend   to $\lim_{q\ra\infty}c_1^2/c_2=1$, and one may say
that our surfaces show no pathological behaviour with respect to surface geography. 
The determination of the Euler characteristic $c_2=e$ depends on Dolgachev's   formula 
$$
e(X) = e(X_{\bar{\eta}}) e(B) + \sum \left(e(X_a) - e(X_{\bar{\eta}}) +\delta_a\right)
$$
for $l$-adic Euler characteristic for schemes   fibered over curves \cite{Dolgachev 1972}, where $\delta_a$
is \emph{Serre's measure of wild ramification}. Its computation is quite easy in our situation, given the
explicit nature of the Artin--Schreier curves. 

The surface $S$ has a surprisingly simple \emph{projective description}, which is   reminiscent 
to  Kummer's quartic surfaces  $A/\left\{\pm 1\right\}\subset\PP^3$.
Passing   to the normal surface $\bar{S}$ with $q-1$ rational double points of type $A_{q-1}$
obtained by contracting all   terminal chains in the fundamental cycle,
the image $\bar{Z}\subset\bar{S}$ of the fundamental cycle remains Cartier, and defines the embedding:

\begin{maintheorem}
The invertible sheaf $\bar{\shL}=\O_{\bar{S}}(\bar{Z})$ is very ample, has $h^0(\bar{\shL})=4$
and embeds the normal surface $\bar{S}$  as a divisor of degree $q$ in $\PP^3$,
sending the rational double points into a line.
\end{maintheorem}

In fact, a canonical divisor is $K_{\bar{S}}=(q-4)\bar{Z}$. So for $q\geq 5$,
the closed embedding $\Phi_{\bar{Z}}:\bar{S}\ra\PP^3$ can be viewed as
an \emph{$m$-canonical map}, for the fractional value $m=1/(q-4)$.  In the simplest case $q=5$, this  is in line
with classification results of Horikawa on minimal surfaces of general type with $K^2=5$ and vanishing
irregularity \cite{Horikawa 1975}.
The projective description also allows us to deduce that our surfaces $S$ admit a lifting into characteristic zero,
at least in the category of algebraic spaces. It would be   interesting to determine
the homogeneous polynomial   describing the image $\bar{S}\subset\PP^3$, but we have made no
attempt to do so.

\medskip
The paper is organized as follows: In Section \ref{Artin--Schreier curves} we review
some relevant facts on Artin--Schreier curve, all of them   well-known.
In Section \ref{Products of Artin--Schreier curves} we study the normal surface $(C\times C')/G$, obtained as the  
quotient of the product of two Artin--Schreier curves with respect to the diagonal action.
We find an explicit equation for the singularitiy, and a dimension formula for the
global sections of the dualizing sheaf. Section \ref{Invariants of the singularity} contains an   analysis
of the minimal resolution of singularities $X\ra(C\times C')/G$. Notable results are
formulas for the fundamental cycle and its genus, as well as some bounds on the arithmetic genus.
In Section \ref{Vanishing of Irregularity} we prove that $H^1(X,\O_X)=0$, such that the Picard scheme is reduced and 0-dimensional.
This relies on a general fact about  group actions with fixed points, which seems to be of independent interest, and
is verified with Grothendieck's theory of $G$-equivariant cohomology.
In Section \ref{Place in the Enriques classification} we determine the place of the smooth surface $X$ in the Enriques classification.
Among other things, this depends on the geometry  of the fibration $X\ra\PP^1$ induced
from the projections on $C\times C'$.
Section \ref{Canonical model and canonical map} contains our analysis of projective models for the surfaces.
Finally, in Section \ref{Numerical invariants and geography} we take up questions from surface geography  and compute   Chern invariants.
This mainly relies on Dolgachev's   formula for $l$-adic Euler characteristics for fibered
schemes.

\begin{acknowledgement}
The first author would like to thank the Mathematisches Institut of the  
Heinrich--Heine--Universit\"{a}t D\"usseldorf, where this work has begun,  for its  warm hospitality.
Research of the first author was partially supported by Grand-in-Aid 
for Scientific Research (C)  20540044, The Ministry of Education, 
Culture, Sports, Science and Technology. We thank the referee for bringing
to our attention   Lorenzini's preprints \cite{Lorenzini 2006}, \cite{Lorenzini 2011a}, \cite{Lorenzini 2011b}.
\end{acknowledgement}

\section{Artin--Schreier curves}
\mylabel{Artin--Schreier curves}

Let $p>0$ be a prime number
and $k$ be an algebraically closed ground field of characteristic $p$.
Consider  \emph{Artin--Schreier curves}  of the form
$$
C:\quad f(y)=x^q-x,
$$
where the left   side of the defining equation is 
a monic polynomial $f(y)=y^{q-1}+\mu_2y^{q-2}+\ldots+\mu_q$ of degree $q-1$ with coefficients from the ground field $k$.
In other words, $C\subset\PP^2$ is defined
by the homogeneous equation
\begin{equation}
\label{homogeneous equation}
Y^{q-1}Z + \mu_2Y^{q-2}Z^2+\ldots+\mu_{q}Z^q = X^q - XZ^{q-1}
\end{equation}
of degree $\deg(C)=q$ inside the projective plane $\PP^2=\Proj(k[X,Y,Z])$. 
Homogeneous and inhomogeneous coordinates are related by $x=X/Z$ and $y=Y/Z$. The curve $C$ is
smooth, with numerical invariants
\begin{equation}
\label{curve genus}
h^0(\O_C)=1 \quadand h^1(\O_C)=(q-1)(q-2)/2\quadand \deg(K_C)=q(q-3).
\end{equation}

Now let $G=\FF_q\subset k$ be the additive group of all scalars 
$\lambda$ satisfying $\lambda^q=\lambda$, 
viewed as an elementary abelian $p$-group.
We may also regard it as a subgroup
$$
G=\left\{
\begin{pmatrix} 
1\\
\lambda & 1
\end{pmatrix}
\mid \lambda\in\FF_q\right\}
\subset\GL(2,k).
$$
Thus the elements $\lambda\in G$ act  on $\PP^2$ via
$X\mapsto X+\lambda Z,\quad Y\mapsto Y,\quad Z\mapsto Z$.
This action leaves the homogeneous equation (\ref{homogeneous equation})
invariant, whence induces an action on $C$.
This action is free, except for a single fixed point $a=(0:1:0)\in C$.
Dehomogenizing  in another way by setting $u=X/Y$, $w=Z/Y$, 
we see that an open neighborhood of the fixed point  is the spectrum of the
coordinate ring
$$
k[u,w]/(u^q-uw^{q-1}-P(w)),
$$
where $P(w)=f(1/w)w^q=w+\mu_2w^2+\ldots+\mu_qw^q$, 
and the group elements $\lambda\in G$ act via $u\mapsto u+\lambda w$, $w\mapsto w$.
Since 
$$
\frac{\partial}{\partial w}(u^q-uw^{q-1}-P(w)) = uw^{q-2}-(1+2\mu_2w+\ldots+ (-\mu_{q-1}w^{q-2}))
$$
becomes a unit in the local ring $\O_{C,a}$, there is
a unique way to write the indeterminate $w$  as a formal power series $w(u)=\sum\alpha_iu^i$ in the variable $u$ so that
$u^q-uw(u)^{q-1}-P(w(u))=0$ (\cite{A 4-7}, \S 4, No.\ 7, Corollary to Proposition 10). 
Using the latter condition, one easily infers that the initial coefficients are
\begin{equation}
\label{formal substitution}
\alpha_0=\ldots=\alpha_{q-1}=0\quadand \alpha_q=1.
\end{equation}
The upshot is that the inclusion
$$
k[[u]]\subset k[[u,w]]/(u^q-uw^{q-1}-P(w))=\O_{C,a}^\wedge
$$
is bijective, and the group elements $\lambda\in G$ act on the formal completion $k[[u]]=\O_{C,a}^\wedge$ via
\begin{equation}
\label{i invariant}
u\longmapsto u + \lambda u^q + \text{higher order terms}.
\end{equation}
From this we infer that the filtration given by the \emph{higher ramification subgroups}
$G=G_0\supset G_1\supset G_2\supset\ldots$
takes the simple form
\begin{equation}
\label{ramification groups}
G_i=
\begin{cases}
G & \text{if $i\leq q-1$};\\
0 & \text{if $i\geq q$}.
\end{cases}
\end{equation}
Recall that $G_i\subset G$ is defined as the \emph{decomposition group} of the $i$-th infinitesimal
neighborhood of the closed point, that is, the subgroup of those
$\sigma\in G$ with the property $\sigma(u)-u\in \maxid^{i+1}_a$.
The corresponding   function $i_G $ on $G$ is given by
$$
i_G(\sigma)=
\begin{cases}
q & \text{if $\sigma\neq 0$};\\
\infty & \text{if $\sigma =0$}.
\end{cases}
$$
We refer to Serre's monograph \cite{Serre 1979} for the theory of higher ramification groups.
These groups will play a crucial  role in Section \ref{Numerical invariants and geography} in the 
determination of Euler characteristics.
The Hurwitz Formula  for the quotient map $C\ra C/G$ of degree $q$,
in the form of \cite{Serre 1979}, Chapter VI, Proposition 7,
gives
$$
2-(q-1)(q-2) = 2-2g_C = q(2-2g_{C/G}) - a_G(0) = q(2-2g_{C/G}) - (q-1)q,
$$
where $g_C$, $g_{C/G}$ denotes genus, and $a_G$ is the character of the Artin representation attached to the fixed point,
which by definition has $a_G(0)=\sum_{\sigma\neq 0}i_G(\sigma)$.
It follows that  $g_{C/G}=0$, whence  $C/G=\PP^1$. In light of the
defining equations, this was of course clear from the very beginning: the quotient map $C\ra\PP^1$ is a classical Artin--Schreier
covering of the projective line. Using that the quotient map is \'etale away from and totally ramified at the fixed point,
one deduces:

\begin{proposition}
\mylabel{canonical curve}
A canonical divisor for the curve $C$ is given by
$K_C=q(q-3)a$, where $a\in C$ is the fixed point for the $G$-action.
\end{proposition}

It is not difficult to give an explicit basis for the vector space of global $1$-forms on $C$:

\begin{proposition}
\mylabel{rational differentials}
The rational differentials $x^iy^jdy$, $0\leq i+j\leq q-3$ are everywhere defined 
and constitute a basis for $H^0(C,\Omega^1_C)$.
\end{proposition}

\proof
Consider the two   coordinate rings
$$
R=k[x,y]/(x^q-x-y^qP(1/y))\quadand R'=k[u,w]/(u^q-uw^{q-1}-P(w))
$$ 
for our curve $C$. The relation $dx=(...)dy$ reveals that
$\Omega^1_R$ is freely  generated by $dy$. 
Similarly, the relation $0=(w^{q-2}u -P'(w))dw -w^{q-1}du$ and $P'(0)=1$
shows that $\Omega^1_{R'}$ is freely generated by $du$, locally at the point $u=w=0$.
Given a polynomial $f(x,y)$,
we thus express the differential $f(x,y)dy$ in terms of $u,w$, using
$y=Y/Z=1/w$ and $x=X/Z=u/w$:
$$
f(x,y)dy = -f(u/w,1/w)w^{-2}dw = \frac{1}{P'(w)-w^{q-2}u} w^{q-3}f(u/w,1/w)du
$$
Hence the rational differential $f(x,y)dy$ is everywhere defined if  
$w^{q-3}f(u/w,1/w)$, which lies in the field of fractions for  $R'$, actually lies in $R'$. 
This indeed holds for the monomials $f(x,y)=x^iy^j$, provided $0\leq i+j\leq q-3$,
by (\ref{formal substitution}).
The resulting $(q-1)(q-2)/2$ elements $x^iy^jdy\in H^0(C,\Omega^1_C)$ are clearly linearly independent over $k$.
They must constitute a basis, because the vector space in question is of dimension $(q-1)(q-2)/2$.
\qed

\medskip
The induced action of the $\lambda\in G$ on $H^0(C,\Omega^1_C)$ is  
$x^iy^jdy\mapsto (x+\lambda)^iy^jdy$.
Whence $H^0(C,\Omega^1_C)$ is the direct sum of the $G$-invariant subspaces
$$
V_j\subset H^0(C,\Omega^1_C),\quad 0\leq j\leq q-3
$$
that are generated by $x^iy^jdy$, $0\leq i\leq q-3-j$. Obviously,   $\dim(V_j)=q-2-j$.

\begin{proposition}
\mylabel{indecomposable subspaces}
The  $G$-invariant subspaces $V_j\subset H^0(C,\Omega^1_C)$ are indecomposable as $G$-modules, 
and the fixed spaces $V_j^G\subset V_j$ are $1$-dimensional.
\end{proposition}

\proof
We first check that the fixed space $V_j^G\subset V_j$ are $1$-dimensional.
Clearly, $y^jdy$ is invariant. Seeking a contradiction, we suppose there is a monic polynomial 
$f(x)$ with $q-3-j\geq\deg(f)>1$
so that the differential $f(x)y^jdy$ is invariant. Factoring $f(x)=\prod(x-\omega_i)$, we see
that the set of roots $\left\{\omega_1,\ldots,\omega_d\right\}$ is invariant
under the substitution $\omega\mapsto\omega+\lambda$, $\lambda\in G$.
Hence $\deg(f)\geq q$, contradiction.

Suppose $V_j$ is decomposable, such that we have decomposition $V_j=V_j'\oplus V_j''$
into nonzero $G$-invariant subspaces. Since $G$ is commutative, there is a basis of $V_j'$
in which all $\lambda\in G$ act via lower triangular matrices (see \cite{A 4-7}, Chapter VII, \S5, No.\ 9 Proposition 19).
The last member $x$ of such a basis  is then a common eigenvector for all $\lambda\in G$.
Since each $\lambda\in G$ has order $p$, all eigenvalues are $\epsilon=1$, whence
$x$ is $G$-fixed. The same applies to $V_j''$, giving a contradiction
to $\dim(V_j^G)=1$.
\qed

\medskip
Note   this implies  that the fixed space $H^0(C,\Omega^1_C)^G$ is of dimension $q-2$.

\section{Products of Artin--Schreier curves}
\label{Products of Artin--Schreier curves}

\medskip
Now choose a second Artin--Schreier curve $C'$ of the form discussed in the previous section,
and consider the product $C\times C'$, endowed with the diagonal $G$-action.
In this section we start to study the quotient
$(C\times C')/G$, which is a normal surface
whose singular locus consists of one point $s\in(C\times C')/G$, the image of the fixed point $(a,a')$.
The projections $\pr_1:C\times C'\ra C$ and $\pr_2:C\times C'\ra C'$ induce fibrations
$$
\varphi_1:(C\times C')/G\ra C/G=\PP^1\quadand
\varphi_2:(C\times C')/G\ra C'/G=\PP^1,
$$ 
respectively. Choose   coordinates on the copies of projective lines so that
the fixed points $a\in C$ and $a'\in C'$ map to the origin $0\in\PP^1$,
and consider the fibers $\varphi_1^{-1}(0),\varphi_2^{-1}(0)\subset (C\times C')/G$.
Recall that the \emph{multiplicity of a fiber} is the greatest common divisor
for the multiplicities of its integral components.

\begin{proposition}
\mylabel{canonical quotient}
The fibers $\varphi_i^{-1}(0)\subset (C\times C')/G$ have multiplicity $q$,
and a canonical divisor is given by $K_{(C\times C')/G}=(q-3)(\varphi_1^{-1}(0)+\varphi_2^{-1}(0))$.
Its selfintersection number is $K_{(C\times C')/G}^2=2q(q-3)^2$.
\end{proposition}

\proof
The fiber $\varphi_i^{-1}(0)$ is clearly irreducible.
According to (\ref{i invariant}), the $G$-action on the preimage  of $0\in\PP^1=C/G$ in $C$ is trivial. 
Whence outside the singularity, the fiber $\varphi_i^{-1}(0)$ is the quotient of
$C\otimes_k \O_{C,a}/\maxid^q$, where the action on the right factor is trivial,
which is an Artin ring of length $q$. It follows that the fiber has multiplicity $q$.

As to the canonical divisor, we have $K_{C\times C'}=q(q-3)(\pr_1^{-1}(a) +\pr_2^{-1}(a'))$
by Proposition \ref{canonical curve}. This is the preimage of $(q-3)(\varphi_1^{-1}(0)+\varphi_2^{-1}(0))$.
Using that the projection $C\times C'\ra(C\times C')/G$ is \'etale in codimension one,
together with  \cite{Peskin 1984}, Theorem 2.7, we infer that $(q-3)(\varphi_1^{-1}(0)+\varphi_2^{-1}(0))$
is indeed a canonical divisor. Its selfintersection number is
$$
K_{(C\times C')/G}^2=\frac{1}{q}K_{C\times C'}^2 = 2q(q-3)^2,
$$
by the Projection Formula.
\qed

\medskip
The canonical divisor $K_{(C\times C')/G}=(q-3)(\varphi_1^{-1}(0)+\varphi_2^{-1}(0))$ is clearly
Cartier, such that the normal surface $(C\times C')/G$ is Gorenstein.
Let us now examine its singularity. Throughout the paper, it will be crucial to understand the local invariants
of this singularity, in order to determine global invariants for smooth models  of $(C\times C')/G$.
Let $u,u',w,w'$ be four variables and set
$$
A=k[[u,u']]=k[[u,w,u',w']]/(u^q-uw^{q-1}-P(w),u'^q-u'w'^{q-1}-Q(w')),
$$
which is the complete local ring at the fixed point $(a,a')\in C\times C'$.
Here
$$
P(w)=w+\mu_2w^2 +\ldots+\mu_qw^q \quadand Q(w')=w'+\mu'_2w'^2 +\ldots+\mu'_qw'^q
$$
are  polynomials stemming from the left hand side of the Artin-Schreier equations, as discussed
in Section \ref{Artin--Schreier curves}. The group elements $\lambda\in G$ act via
$$
u\longmapsto u+\lambda w,\quad u'\longmapsto u'+\lambda w',\quad w\longmapsto w,\quad w'\longmapsto w'.
$$
Clearly, the elements $w,w',wu'-w'u\in A$ are invariant and satisfy the relation
$$
(wu'-w'u)^q=(ww')^{q-1}(wu'-w'u)+w^qQ(w')-w'^qP(w).
$$
Therefore, we obtain a homomorphism of $k$-algebras
\begin{equation}
\label{quotient singularity}
k[[a,b,c]]/(c^q-(ab)^{q-1}c-a^qQ(b)+b^qP(a)) \lra A^G,
\end{equation}
where $a,b,c$ are indeterminates and $a\mapsto w$, $b\mapsto w'$, $c\mapsto wu'-w'u$.

\begin{proposition}
\mylabel{invariant ring}
The preceding homomorphism (\ref{quotient singularity}) is bijective.
\end{proposition}

\proof
Clearly, the local ring  $R=k[[a,b,c]]/(c^q-(ab)^{q-1}c-a^qQ(b)+b^qP(a))$
is 2-dimensional  and Cohen--Macaulay. Computing the jacobian ideal, one
sees that the singular locus consists of the closed point, whence
$R$ is normal. According to Galois Theory, the finite extension $A^G\subset A$
has generically rank $q=\ord(G)$. By the Main Theorem of Zariski (\cite{EGA IIIb}, Corollary 4.4.9), it therefore
suffices to check that $R\subset A$ has generically rank $q$.
Obviously, the extensions $k[[a,b]]\subset R$ and $k[[a,b]]\subset A$ have
rank $q$ and $q^2$, respectively, and the statement follows by transitivity of ranks.
\qed

\medskip
This local description of the singularity enables us to determine the scheme structure for 
the canonical divisor $K_{(C\times C')/G}= (q-3)(\varphi_1^{-1}(0)+\varphi_2^{-1}(0))$, viewed as a reduced 
Weil divisor:

\begin{corollary}
\mylabel{fiber union}
The reduced Weil divisors $\varphi_i^{-1}(0)_\red$ are isomorphic to the projective line,
and the schematic intersection $\varphi_1^{-1}(0)_\red\cap\varphi_2^{-1}(0)_\red$ has length one.
\end{corollary}

\proof
We saw in the proof of Proposition \ref{canonical quotient} that $\varphi_i^{-1}(0)_\red$ is isomorphic
to $C/G$, at least outside the singular point on the ambient surface.
However, at this singularity, $\varphi_i^{-1}(0)$ is formally isomorphic to the spectrum of 
$k[[a,b,c]]/(c^q, a)$, which follows from Proposition \ref{invariant ring}. Hence its reduction is regular.
Similarly, the union $\varphi_1^{-1}(0)_\red\cup\varphi_2^{-1}(0)_\red$ is formally isomorphic
to the spectrum of $k[[a,b,c]]/(c,ab)$, and the result follows.
\qed

\medskip
We would like to compute the dimension of the vector space of global  sections for the dualizing
sheaf $\omega_{(C\times C')/G}$, but are only able to do so directly in the special case $q=p$:

\begin{proposition}
\mylabel{invariant differentials}
For $q=p$, we have $h^0(\omega_{(C\times C')/G})=(2p^3-9p^2+13p-6)/6=(2p-3)(p-2)(p-1)/6$.
\end{proposition}

\proof
We have a commutative diagram
$$
\begin{CD}
H^0(U,\omega_{C\times C'})^G         @<<< H^0(\omega_{C\times C'})^G\\
@AAA @AAA\\
H^0(V,\omega_{(C\times C')/G}) @<<< H^0(\omega_{(C\times C')/G})
\end{CD}
$$
where $U\subset C\times C'$ is the locus where $G$ acts freely,
$V\subset(C\times C')/G$ is the smooth locus, and $U\ra V$ is the induced finite \'etale
Galois covering with Galois group $G$.
The latter ensures that the vertical map on the left is bijective.
Since our dualizing sheaves are invertible, the horizontal maps are bijective 
(compare, for example, \cite{Hartshorne 1994}, Proposition 1.11). It follows that that the canonical map
$H^0(\omega_{(C\times C')/G})\ra H^0(\omega_{C\times C'})^G$ is bijective.

Recall that $H^0(C,\omega_C)=\bigoplus_{j=0}^{p-3}V_j$, where $V_j$ 
are indecomposable $G$-submodule of dimension $d=p-2-j$.
Since $G=\FF_p$, the action of the generator $1\in\FF_p$ is given, in a suitable
basis, by the Jordan matrix
$$
J_d(1)=\begin{pmatrix}
1 \\
1 & 1 \\
  & \ddots & \ddots\\
  &         & 1 & 1
\end{pmatrix}
\in\Mat(d,k),
$$
which determines the $G$-module up to isomorphism.
Computing the dimension of the $G$-fixed part for $H^0(\omega_{C\times C'})=H^0(C,\omega_C)\otimes H^0(C',\omega_{C'})$ 
thus reduces to extracting the number of  blocks in the Jordan normal form of tensor products of certain Jordan matrices. Write
$$
J_d(1)\otimes J_{d'}(1)=\bigoplus_r  J_r(1)^{\oplus\lambda_r}
$$
for certain multiplicities $\lambda_r\geq 0$. Note that it is a notorious unsolved problem in linear algebra to 
find the Jordan decomposition of tensor products of nilpotent Jordan matrices in positive characteristics.
However, it is well-known that $\sum_r\lambda_r=\min(d,d')$ (see, for example, \cite{Schroeer 2010}, Proposition 3.2).
Hence the $G$-fixed part of such a tensor product representation is   of dimension $\min(d,d')$, and we get
$$
h^0(\omega_{(C\times C')/G})=\sum_{d,d'=1}^{p-2} \min(d,d').
$$
Summing over $l=\min(d,d')$ rather than $(d,d')$, we rewrite the latter sum as
$$
\sum_{l=1}^{p-2}(2(p-2-l)+1)l = \sum_{l=1}^{p-2}(-2l^2+(2p-3)l).
$$
The statement follows by applying the formulas $1+2+\ldots+n=n(n+1)/2$ and $1^2+2^2+\ldots +n^2=n(n+1)(2n+1)/6$.
\qed

\begin{remark}
In order to make similar computations in the general case $q=p^s$, one would need to understand the modular representations
of the elementary abelian group $\FF_q=(\ZZ/p\ZZ)^s$. Such representations can be expressed in terms of $s$ commuting nilpotent matrices,
or equivalently via modules of finite length over the polynomial ring in $s$ indeterminates over $\FF_p$.
For this situation, little seems to be known about multiplicities in tensor products.
\end{remark}

\section{Invariants of the singularity}
\label{Invariants of the singularity}

In this section we study in more detail the  2-dimensional ring
$$
R=k[a,b,c]/(c^q-(ab)^{q-1}c-a^qQ(b)+b^qP(a)),
$$
whose formal completion gives the singularity of the surface
$(C\times C')/G$.
As in the proof of Proposition \ref{invariant ring}, one sees that $R$ is  normal,
and the maximal ideal $\maxid=(a,b,c)$ corresponds 
the unique singularity   $s\in\Spec(R)$.
According to \cite{Lorenzini 2006}, Theorem 2.5, the  exceptional divisor on the minimal resolution of singularities
consists of projective lines, and has as dual graph a tree. To understand its structure, we first consider a partial resolution of singularities:
Let 
$$
f:Y\lra\Spec(R)
$$ 
be the blowing-up of the maximal ideal $\maxid\subset R$, and denote by $E=f^{-1}(s)$  the exceptional divisor. This is the Cartier divisor $E\subset Y$
with ideal $\O_Y(1)\subset\O_Y$.

\begin{proposition}
\mylabel{first blowing-up}
The  surface $Y$ is normal, and $\Sing(Y)=\left\{s_0,\ldots,s_{q}\right\}$ consists of $q+1$
closed points, whose local rings are rational double points of type $A_{q-1}$. 
We have $E_\red=\PP^1$ and $qE_\red=E$. Furthermore, the   inclusions
$iE_\red\subset (i+1)E_{\red}$ are infinitesimal extensions by the sheaf $\O_{\PP^1}(-i)$, $0\leq i\leq q-1$.
\end{proposition}

\proof
The blowing-up $Y$ is covered by three affine charts, the $a$-chart, the $b$-chart, and the $c$-chart.
The coordinate ring of the $a$-chart is generated by $a,b/a,c/a$, subject to the relation
$$
(c/a)^q-a^{q-1}(b/a)^{q-1}c/a - Q(a\cdot b/a)  + (b/a)^qP(a)=0.
$$ 
Computing partial derivatives, one sees that the singular locus
is given by $(c/a) = a = 0$  and $(b/a)^q=(b/a)$, and the latter means $b/a=i$, $i\in\FF_q$. Writing
the left hand side of the preceding equation in the form
$$
\label{rdp equation}
(c/a)^q - a\Psi(c/a,b/a,a),
$$
one easily computes that
$c/a,a,\Psi$ form a regular system of parameters in the formal completion $k[[c/a,b/a-i,a]]$.
The upshot is that  singularity is a rational double point of type $A_{q-1}$.
The situation on the $b$-chart is symmetric, whereas the $c$-chart turns out
to be disjoint from the exceptional divisor.

On the $a$-chart, the  exceptional divisor is given by $a=0$,
whence its coordinate ring is $k[b/a,c/a]/(c/a)^q$. Its reduction is defined by $c/a=0$.
We infer that $E_\red$ is isomorphic to $\PP^1=\Proj k[a/b,b/a]$. The ideal of $iE_\red\subset (i+1)E_{\red}$
is generated by $(c/a)^i$ and $(c/b)^i$ on the $a$- and  $b$-chart, respectively. The statement about
the   infinitesimal extensions follows.
\qed

\medskip
Now let $g:X\ra Y$ be the minimal resolution of singularities of $Y$, the latter being a partial resolution of $\Spec(R)$.  
For each singular point $s_i\in Y$, $0\leq i\leq q$, write
$g^{-1}(s_i)=A_{i,1}\cup\ldots\cup A_{i,q-1}$ as a chain of projective lines $A_{i,j}$, such that  $A_{i,j}\cdot A_{i,j+1}=1$.
Let $A_0\subset X$ be the strict transform of $E_\red\subset Y$, which is another projective line.

According to Lorenzini \cite{Lorenzini 2006}, Theorem 2.5, the dual graph for the exceptional divisor on the resolution 
of a wild quotient singularity is necessarily a tree. To our knowledge, no examples are known where the tree is not star-shaped.
In \cite{Lorenzini 2011a}, Question 1.1, Lorenzini askes whether trees with more that one node are possible.
In our situation, Equation (\ref{rdp equation}) immediately shows that 
$A_0$ intersects each chain of rational curves $g^{-1}(s_i)$ in a terminal component of the chain, 
say $A_{i,1}$. Thus the reduced exceptional divisor    
$$
A_0 + \sum_{i=0}^q \sum_{j=1}^{q-1} A_{i,j}
$$
is a strictly normal crossing made out of projective lines, and its dual graph looks like this:

\vspace{2em}
\centerline{\includegraphics{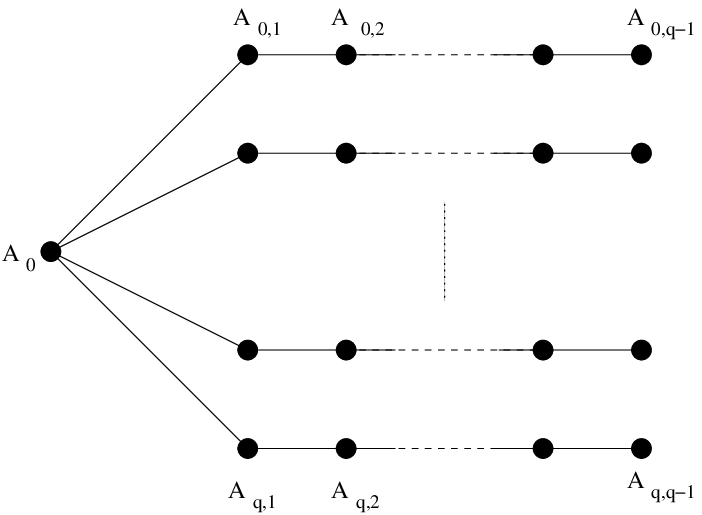}} 
\vspace{1em}
\centerline{Figure \stepcounter{figure}\arabic{figure}: Dual graph for exceptional divisor on $X$.}
\vspace{1em}

\begin{proposition}
\mylabel{selfintersection numbers}
We have $E^2=-q$, and $E\cdot E_\red =-1$,  and $A_0^2= -q$. In particular,
the composite map $X\ra\Spec(R)$ is the minimal resolution of singularities.
\end{proposition}

\proof
As   in the proof of Proposition \ref{first blowing-up},  we have
$$
E_\red =\Spec k[a/b] \cup \Spec k[b/a],
$$
and the invertible sheaf $\O_Y(1)|_{E_\red}$ is given by $a=a/b\cdot b$ on the two charts.
Viewing $a/b\in k[a/b,b/a]^\times$ as a cocycle for $\O_Y(1)|_{E_\red}$, one immediately verifies
that this invertible sheaf has degree $1$. Thus 
$$
E\cdot E_\red =\deg\O_{E_\red}(E) = \deg\O_{E_\red}(-1)=-1.
$$
Since $E=qE_\red$,  the selfintersection number $E^2=-q$ also follows.

Now consider the first infinitesimal neighborhood $E_\red\subset 2E_\red\subset E$.
In light of Proposition \ref{first blowing-up},
this is an infinitesimal extension of $E_\red=\PP^1$ by the invertible sheaf $\shL=\O_{\PP^1}(-1)$.
Let $Y'\ra Y$ be the blowing-up of the singular points $s_0,\ldots, s_q\in Y$.
Then the exceptional divisors are disjoint unions of pairs of rational curves, each pair intersecting transversely at one
point whose local ring on $Y'$ is a rational double point of type $A_{q-3}$.
Moreover, the strict transform of $E_\red$ lies in the smooth locus of $Y'$.
According to the theory of \emph{ribbons} developed by Bayer and Eisenbud, the strict transform of $2E_\red$ on $Y'$,
that is, the blowing-up of the scheme $2E_\red$ with respect to the centers $s_0,\ldots,s_q\in 2E_\red$,
is an infinitesimal extension of $\PP^1$ by the invertible sheaf $\shL(\sum a_i)=\O_{\PP^1}(q)$, compare \cite{Bayer; Eisenbud 1995}, Theorem 1.9.
It follows that $A_0^2=\deg(\O_{A_0}(A_0))=-q$.  In particular, the exceptional divisor for the
resolution of singularities $X\ra\Spec(R)$ contains no $(-1)$-curve, hence is minimal.
\qed

\medskip
We next compute the \emph{fundamental cycle} $Z\subset X$ for the resolution of singularities
$h:X\ra\Spec(R)$, a notion introduced by M.\ Artin \cite{Artin 1966}.
By definition, this is the smallest   effective cycle $Z$ whose support equals the exceptional
divisor and has intersection number $\leq 0$ on each irreducible component of the
exceptional divisor. One way to compute  the fundamental cycle is with a \emph{computation sequence}
$Z_\red=Z_0\subset Z_1\subset\ldots\subset Z_r=Z$, where  in each step $Z_{s+1}-Z_s$ is
an integral component of the exceptional divisor with  $Z_s\cdot (Z_{s+1}-Z_s)>0$.
Using the exact sequence
\begin{equation}
\label{sucsessive extension}
0\lra \O_{Z_{s+1}-Z_s}(-Z_s)\lra\O_{Z_{s+1}}\lra\O_{Z_s}\lra 0,
\end{equation}
one inductively infers that $h^0(\O_{Z_{s+1}})=h^0(\O_{Z_s})$, in particular $h^0(\O_Z)=1$ .
Consequently the schematic image $h(Z)\subset\Spec(R)$ is nothing but the reduced  singular point $s\in\Spec(R)$.
Indeed, one should view the fundamental cycle as an  approximation to the schematic fiber $h^{-1}(s)\subset X$,
the latter usually containing embedded components.

\begin{proposition}
\mylabel{fundamental cycle}
The fundamental cycle is given by the formula 
$$
Z= qA_0+\sum_{i=0}^q \sum_{j=1}^{q-1} (q-j)A_{ij},
$$
and its selfintersection number is $Z^2=-q$.
\end{proposition}

\proof
Let $Z'$ be the cycle on the right hand side. One easily computes the intersection numbers
\begin{equation}
\label{intersection fundamental}
\begin{split}
Z'\cdot A_0 &= -q^2+(q+1)(q-1) = -1,\\
Z'\cdot A_{ij}&= (q-j+1) -2(q-j)+(q-j-1)=0,
\end{split}
\end{equation}
whence the fundamental cycle $Z$ is contained in   $Z'$, by the
minimality property of fundamental cycles.
Seeking a contradiction, we suppose $Z\subsetneqq Z'$.
The effective cycle $Z'-Z$ has intersection numbers
\begin{align*}
(Z'-Z)\cdot A_0 &= -1 - (Z\cdot A_0),  \\
(Z'-Z)\cdot A_{ij}         &= -Z\cdot A_{ij}\geq 0.
\end{align*}
Since nonzero effective exceptional cycles are not nef on all exceptional curves
(see \cite{Artin 1966}, proof of Proposition 2),
there is no possibility but $Z\cdot A_0=0$. Now write
$Z=\lambda_0A_0 + \sum_{i=0}^q\sum_{j=1}^{q-1}\lambda_j A_{ij}$
with coefficients $1\leq \lambda_0\leq q$ and $1\leq\lambda_j\leq q-j$.
Note that the coefficients $\lambda_j$ do not depend on $i$, due to the obvious symmetry of
the dual graph in Figure 1. We have
$$
0=Z\cdot A_0 = -q\lambda_0 + (q+1)\lambda_1,
$$
whence $q\mid \lambda_1$, contradicting $1\leq \lambda_1\leq q-1$.
\qed


\medskip
Recall that  the \emph{canonical cycle} $K_h=K_{X/R}$ is the divisor
supported by the exceptional divisor that satisfies the equations $K_{X/R} \cdot C+C^2 = \deg(K_C)$,
where $C$ runs through the integral exceptional divisors. 

\begin{corollary}
\mylabel{canonical cycle}
The canonical cycle is given by $K_h= -(q-2)Z$, with selfintersection number 
$K_h^2=-q(q-2)^2$.
\end{corollary}

\proof
Using that $A_0$ and $A_{i,j}$ are copies of the projective lines
with selfintersection  numbers $A_0^2=-q$, $A_{i,j}^2=-2$, one deduces the result from
the intersection numbers
(\ref{intersection fundamental}). The selfintersection then follows from Proposition \ref{fundamental cycle}.
\qed

\medskip
We are now in position to compute the arithmetic genus of the fundamental cycle
$$
p_f = h^1(\O_Z)= 1-\chi(\O_Z),
$$
which is also called the \emph{fundamental genus} of the singularity, confer \cite{Tomaru 1995}.

\begin{corollary}
\mylabel{fundamental genus}
The  fundamental genus of the singularity $R$ is given by the formula $p_f=(q-1)(q-2)/2$.
\end{corollary}

\proof
Riemann--Roch yields $-2\chi(\O_Z)=\deg(K_Z) = Z^2 + K_h\cdot Z=(3-q)Z^2$.
We have $Z^2= -q$, and finally obtain
$$
p_f = 1-\chi(\O_Z)= 1-(3-q)q/2=(q-1)(q-2)/2,
$$
as claimed.
\qed

\medskip
From this we deduce:

\begin{corollary}
\mylabel{rdp elliptic}
The singularity $R$ is a rational double point if and only if $q=2$,
and minimally elliptic if and only if $q=3$.
\end{corollary}

\proof
According to \cite{Artin 1966}, Theorem 3,  the singularity is rational if and only if $h^1(\O_Z)=0$.
In light of Corollary \ref{fundamental genus}, this happens precisely when $q=2$.
By Laufer \cite{Laufer 1977}, Theorem 3.4, the condition $K_h=-Z$ is one of several
equivalent defining property of \emph{minimally elliptic singularities}. By Corollary \ref{canonical cycle},
this happens if and only if $q=3$.
\qed

\begin{remark}
For  $q=2$, this singularity is actually a rational double point of type $D_4^1$,
according to Artin's list \cite{Artin 1977}. Indeed, the equation $c^q-(ab)^{q-1}c-a^qb+b^qa=0$ is a special case of
Artin's normal form for wild $\ZZ/2\ZZ$-quotient singularities in dimension two \cite{Artin 1975},
compare also \cite{Schroeer 2009}. For $q=3$, the minimally elliptic singularity appears in 
Laufer's classification  (\cite{Laufer 1977}, Table 3 on page 1294)
under the designation $A_{1,\star,0}+A_{1,\star,0}+A_{1,\star,0}+A_{1,\star,0}$.
\end{remark}

\begin{remark}
Shioda \cite{Shioda 1974} and Katsura \cite{Katsura 1978} obtained rather similar results for the action of the sign involution   
on abelian surfaces in characteristic $p=2$.
\end{remark}

We next want to compute the  \emph{geometric genus}
$$
p_g= \length R^1h_*(\O_X) = h^1(\O_{nZ}),\quad n\gg 0
$$
of the singularity. Except for rational double points and minimally elliptic singularities,
this invariant is difficult to compute. We have at least some bounds.
It will turn out later that in the special case $q=p$,
these bounds are actually equalities (Corollary \ref{local geometric genus}).

\begin{proposition}
\mylabel{geometric genus}
The geometric genus of the singularity $R$ satisfies the inequalities $p_g\leq q(q-1)(q-2)/6$.
For $q\geq 5$, we moreover have $h^1(\O_{2Z})>h^1(\O_Z)$, in particular  $p_g>p_f$.
\end{proposition}

\proof
It is   convenient to work on the partial resolution $Y$ rather than on the full resolution $X$.
The Leray--Serre spectral sequence gives an exact sequence
$$
0\lra R^1f_*(\O_Y) \lra R^1h_*(\O_X) \lra f_*(R^1g_*(\O_X)).
$$
The term on the right vanishes, since $Y$ has only rational singularities,
whence $p_g=\length R^1f_*(\O_Y)$.
Let $D=-K_f=(q-2)E$ be the anticanonical cycle, where $E=f^{-1}(s)$. It follows with the Grauert--Riemenschneider Vanishing Theorem
(see \cite{Giraud 1982}, Theorem 1.5) that $R^1f_*(\O_Y(-D))=R^1h_*(g^*\O_Y(-D))=0$, whence the
canonical surjection 
$$
H^0(R^1f_*(\O_Y))\lra H^1(D,\O_D)
$$
is bijective. Moreover, Riemann--Roch gives $\chi(\O_D)= (D+K_f)\cdot D=0$, such
that $p_g=h^1(\O_D)=h^0(\O_D)$.

According to Proposition  \ref{first blowing-up} and Proposition \ref{selfintersection numbers},
we have 
$$
E_\red=\PP^1\quadand E=qE_\red\quadand E\cdot E_\red =-1.
$$
Consider the integral Weil divisors $(i/q)E=iE_\red$.
We deduce from Proposition \ref{first blowing-up} that the kernel $\shK_i$
in the exact sequence
$$
0\lra \shK_i \lra \O_{((i+1)/q)E} \lra \O_{(i/q)E}\lra 0
$$
is an invertible sheaf on $E_\red=\PP^1$ of degree
$$
\deg(\shK_i) = \lfloor i/q \rfloor -q \left\{i/q\right\} 
$$
where $\lfloor i/q\rfloor$ and $\left\{i/q\right\}$ denotes   integral and fractional parts, respectively.
Let us tabulate the kernels $\shK_i$ for $0\leq i< q(q-2)$
in a matrix of size $(q-2)\times q$:
\begin{equation}
\label{kernel array}
\begin{array}{ccccccc}
\O_{\PP^1}(0)   & \O_{\PP^1}(-1)  & \O_{\PP^1}(-2) & \quad&\ldots         &                & \O_{\PP^1}(1-q)\\
\O_{\PP^1}(1)   & \O_{\PP^1}(0)   & \O_{\PP^1}(-1) & &\ldots         &                & \O_{\PP^1}(2-q)\\
\vdots          & \vdots          &                & & \ddots        &                & \vdots\\
\O_{\PP^1}(q-3) & \O_{\PP^1}(q-4) & \ldots         & &               & \O_{\PP^1}(-1) & \O_{\PP^1}(-2)
\end{array}
\end{equation}
Only kernels of degree $\geq 0$ may contribute to $h^0(\O_D)$, and the total possible contribution is
$$
\begin{aligned}
 & (q-2)\cdot h^0(\O_{\PP^1}) + (q-3)\cdot h^0(\O_{\PP^1}(1)) + \ldots + 1\cdot h^0(\O_{\PP^1}(q-3))\\
 & = \sum_{i=1}^{q-2} i(q-1-i)\\
 & = (q-1)\sum_{i=1}^{q-2} i - \sum_{i=1}^{q-2} i^2\\
 & = (q-1)(q-1)(q-2)/2 - (q-2)(q-1)(2q-3)/6\\
 & = q(q-1)(q-2)/6.
\end{aligned}
$$
But some  coboundary maps in the exact sequence
\begin{equation}
\label{coboundary maps}
H^0(\O_{((i+1)/q)E})\lra H^0(\O_{(i/q)E})\stackrel{\partial}{\lra} H^1(\shK_i) \lra H^1(\O_{((i+1)/q)E})
\end{equation}
might be nonzero, so we only get an upper bound $p_g\leq q(q-1)(q-2)/6$, rather than an equality.

To obtain a lower bound, we consider the  cycle $2E\subset Y$. The second row in (\ref{kernel array})
reveals that some cycle  $E\subset F\subset 2E$ has
$$
h^0(\O_{F})=h^0(\O_{\PP^1}(0)) + h^0(\O_{\PP^1}(1)) +h^0(\O_{\PP^1}(0)) =4 \quadand h^1(\O_{F})=h^1(\O_E).
$$
Now suppose that we would have $h^1(\O_F)=h^1(\O_{2E})$. Then in each step leading from $F$ to $2E$ the
coboundary map  in (\ref{coboundary maps}) must surjects onto the respective cohomology groups
$$
H^1(\O_{\PP^1}(-2)),\quad H^1(\O_{\PP^1}(-3)),\quad \ldots \quad H^1(\O_{\PP^1}(2-q)).
$$
It follows that 
$$
4=h^0(\O_{F}) >   h^1(\O_{\PP^1}(2-q)) = q-3,
$$
whence $q\leq 5$. This already implies the inequality $p_f<p_g$ for $q>5$.

It remains to rule out the case    $h^1(\O_E)=h^1(\O_{2E})$ and $q=5$.
The preceding paragraph reveals that than $h^0(\O_{2E})=1$, and the short exact sequence
$$
0\lra \O_E(-E)\lra \O_{2E}\lra\O_E\lra 0
$$
yields $\chi(\O_E(-E))=0$. On the other hand, Riemann--Roch gives
$$
\chi(\O_E(-E))=-Z^2+\chi(\O_E) = -q + 1 - (q-1)(q-2)/2=-q(q+1)/2<0,
$$
contradiction.
\qed

\section{Vanishing of Irregularity}
\label{Vanishing of Irregularity}

We continue to study the normal surface $(C\times C')/G$.
Let $f:Y\ra (C\times C')/G$ be the blowing-up of the unique
singularity. We saw in the preceding section  
that the singular locus $\Sing(Y)$ consists of $q+1$ rational double points
of type $A_{q-1}$. Let $g:X\ra Y$ be the minimal resolution of these double points.
We now dispose off the \emph{irregularity} $h^1(\O_X)$. 

\begin{proposition}
\mylabel{irregularity}
The irregularity $h^1(\O_X)$ vanishes.
\end{proposition}

\proof
Since $Y$ contains only rational singularities, it suffices to check   $h^1(\O_Y)=0$.
Let $\tilde{Y}\ra C\times C'$ be the blowing-up of the reduced fixed point $(a,a')$.
We claim that the schematic preimage on $\tilde{Y}$ of the singular point $s\in (C\times C')/G$ is 
a Cartier divisor, such that the universal property of blowing-ups gives a commutative diagram
$$
\begin{CD}
C\times C' & @<<< \tilde{Y}\\
@VVV & @VVV\\
(C\times C')/G & @<<< Y.
\end{CD}
$$
This is a local problem. Using the notation from Section \ref{Products of Artin--Schreier curves},
we have to understand what happens with the ideal
$$
(w,w')= (w,w',wu'-w'u)\subset k[[u,u']] 
$$
on $\tilde{Y}$. But  $(w,w')=(u^q,u'^q)$ according to (\ref{formal substitution}), and it obvious that the latter
ideal becomes invertible upon blowing-up of $(u,u')$.

The $G$-action on $C\times C'$ induces a $G$-action on $\tilde{Y}$, and
the canonical morphism $\tilde{Y}/G\ra Y$ is an isomorphism by the Main Theorem of Zariski.
It follows from  (\ref{formal substitution}) that the $G$-action on the cotangent space $\maxid/\maxid^2$
for the fixed point $(a,a')\in C\times C'$ is trivial, whence the $G$-action
on the exceptional curve $\tilde{E}\subset\tilde{Y}$ is trivial as well.
Now Lemma \ref{injective pullback} below ensures that the induced map $H^1(Y,\O_Y)\ra H^1(\tilde{Y},\O_{\tilde{Y}})$
is injective.

To finish the argument, consider the reduced fiber union $F=\varphi_1^{-1}(0)_\red\cup\varphi_2^{-1}(0)_\red$
for the two projections $\varphi_i:(C\times C')/G\ra\PP^1$. Its preimage on $\tilde{Y}$
contains the strict transform $\tilde{F}\subset\tilde{Y}$ of $C\times\left\{a'\right\}\cup \left\{a\right\}\times C'\subset C\times C'$,
which is isomorphic to a disjoint union $C\amalg C'$, and we have a commutative diagram
$$
\begin{CD}
H^1(\tilde{Y},\O_{\tilde{Y}}) @>>> H^1(\tilde{F},\O_{\tilde{F}})\\
@AAA @AAA\\
H^1(Y,\O_Y) @>>>H^1(F,\O_F).
\end{CD}
$$
The term $H^1(F,\O_F)$ vanishes by Corollary \ref{fiber union}.
Since the  map on the left is injective, it suffices to check that
the the restriction map $H^1(\tilde{Y},\O_{\tilde{Y}})\ra H^1(\tilde{F},\O_{\tilde{F}})$
is injective. Indeed, the maps
$$
H^1(C,\O_C)\oplus H^1(C',\O_{C'})=H^1(\tilde{Y},\O_{\tilde{Y}})\lra H^1(\tilde{F},\O_{\tilde{F}})= H^1(C,\O_C)\oplus H^1(C',\O_{C'})
$$
are all bijective.
\qed

\medskip
In the course of the preceding proof we have used a fact that appears to be of independent interest.
Let us formulate it in a rather general way:
Suppose $X$ is a   scheme over a field $k$,   and   $G$ be a finite group
acting on $X$ so that the quotient $Y=X/G$ exists as a scheme.

\begin{lemma}
\mylabel{injective pullback}
Assumptions as above. Suppose additionally that $k=H^0(X,\O_X)$ and that
there is a rational fixed point $x\in X$. Then the canonical map 
$H^1(Y,\O_Y)\ra H^1(X,\O_X)$ is injective.
\end{lemma}

\proof
The idea is to use $G$-equivariant cohomology $H^r(X,G,\O_X)$, which was introduced in \cite{Grothendieck 1957}.
Consider the two spectral sequences with $E_2$-terms
$$
E_2^{r,s}=H^r(Y,\shH^s(G,\O_X)) \quadand E_2^{r,s}= H^r(G,H^s(X,\O_X))
$$
abutting to $H^{r+s}(X,G,\O_X)$, where $\shH^s(G,\O_X)$ denotes the sheaf of cohomology groups. They give rise to a 
commutative  diagram 
\begin{equation}
\label{spectral sequence}
\begin{xy}
\xymatrix{
        &                                 & 0\ar[d]\\
        &                                 & H^1(Y,\shH^0(G,\O_X))\ar[d]\ar[dr]\\
0\ar[r] & H^1(G,H^0(X,\O_X))\ar[r]\ar[dr] & H^1(X,G,\O_X)\ar[r]\ar[d]          & H^0(G, H^1(X,\O_X))\\
        &                                 & H^0(G,\shH^1(G,\O_X))\\
}
\end{xy}
\end{equation}
with exact row and column, and the composition  given by the upper diagonal arrow
$$
H^1(Y,\O_Y)=H^1(Y,\shH^0(G,\O_X))\lra H^0(G,H^1(X,\O_X))\subset H^1(X,\O_X)
$$ 
is our map in question.
By a diagram chase, it therefore suffices to check that the other composition
$H^1(G,H^0(X,\O_X))\ra H^0(G,\shH^1(G,\O_X))$ is injective.
Now comes in our rational fixed point $x\in X$: Composing further with the restriction map induced
by $\left\{x\right\}\subset X$,
we obtain 
$$
H^1(G,H^0(X,\O_X))\lra H^0(G,\shH^1(G,\kappa(x)))=H^1(G,\kappa(x)).
$$
By assumption, the map $H^0(X,\O_X)\ra\kappa(x)$ is bijective, whence the assertion.
\qed

\begin{remark} 
Let $\nu:X\ra Y$ be the quotient map.
If $X$ is normal, and the order of $G$ is prime to the characteristic of $k$, 
then the existence of a trace map shows that $\O_Y\subset\nu_*(\O_X)$  
is a direct summand, such that $H^r(Y,\O_Y)\ra H^r(X,\O_X)$ is injective
for all $r\geq 0$.
\end{remark}

\begin{remark}
On the other hand, if $Y$ is an Enriques surface in characteristic $p=2$
with $\Pic^\tau_Y=\mu_2$, and $X\ra Y$ is the K3-covering,  such that $G=\pi_1(Y)$ is cyclic of order two
and $Y=X/G$, then $H^1(Y,\O_Y)$ is $1$-dimensional,
whereas $H^1(X,\O_X)$ vanishes. We note in passing that this situation is somewhat typical:
\end{remark}

\begin{lemma}
\mylabel{almost injective}
Suppose $G$ is cyclic, acts freely on $X$, and $k=H^0(X,\O_X)$.
Then the kernel of the canonical map $H^1(Y,\O_Y)\ra H^1(X,\O_X)$ is at most
$1$-dimensional. 
\end{lemma}

\proof
Making a diagram chase in (\ref{spectral sequence}), the dimension of the kernel
is bounded by the dimension of $H^1(G, H^0(X,\O_X))$. Clearly, $H^0(X,\O_X)=k$ is the trivial $G$-module.
Let $n=\ord(G)$. Then $H^1(G,k)$ is isomorphic to the kernel of the multiplication map $n:k\ra k$, whence 
a $k$-vector space of dimension at most one.
\qed

\section{Place in the Enriques classification}
\label{Place in the Enriques classification}

We now study the global geometry of the normal surface $(C\times C')/G$ in more detail.
Recall that $f:Y\ra (C\times C')/G$ is the blowing-up of the singularity. 
We saw in Section \ref{Invariants of the singularity}
that $Y$ is normal  and contains $(q-1)$ rational double points of type $A_{q-1}$.
Let $g:X\ra Y$ be the minimal resolution of these singularities. Then the composite map
$h:X\ra (C\times C')/G$ is the minimal resolution of singularities.
Let $X\ra S$ be the contraction to a minimal model $S$. We display our surfaces and maps
in a commutative diagram:
$$
\begin{xy}
\xymatrix{
C\times C'\ar[dr] &                       &                      & X\ar[dl]^{g}\ar[rr]\ar[dll]_{h}\ar[dd]^{\psi_1,\psi_2} 
 && S \\
                  & (C\times C')/G\ar[drr]_{\varphi_1,\varphi_2} & Y\ar[l]^{ \quad\qquad f}\\
                  &                                              &             &\PP^1
}
\end{xy}
$$
The top row contains the smooth surfaces, the middle row the normal surfaces, and the arrows $\varphi_i$, $\psi_i$
are the maps induced from the two projections $\pr_1:C\times C'\ra C$ and $\pr_2:C\times C'\ra C'$.

The goal of this section is to determine 
the place of $X$, or rather  its minimal model $S$, in the Enriques classification of surfaces, in dependence on the
prime power $q=p^s$. An elementary argument involving only intersection numbers
already gives:

\begin{proposition}
\mylabel{X general type}
If $q\geq 7$, then the  surface $X$ is of general type.
\end{proposition}

\proof
We first compute the number $K_X^2$ on the smooth surface $X$.
Obviously $K_X=K_h + h^*K_{(C\times C')/G}$, such that
$K_X^2=K_h^2 + K_{(C\times C')/G}^2$.
But $K_h^2=-q(q-2)^2$ by Proposition \ref{canonical cycle},
whereas $K_{(C\times C')/G}^2 = 2q(q-3)^2$ according to Proposition \ref{canonical quotient}.
The upshot is
\begin{equation}
\label{selfintersection resolution}
K_X^2 = q(q^2-8q+14) = q((q-4)^2-2).
\end{equation}
Now suppose that $q\geq 7$. Then $K_X^2>0$, and the Theorem of Riemann--Roch gives
$$
\chi(\omega_X^{\otimes t})= (t^2-t)K_X^2/2 +\chi(\O_X).
$$
So for $t\gg 0$   either $tK_X$ or $(1-t)K_X$ is effective.
The latter is impossible, because the canonical divisor on $X$ maps
to the canonical divisor on $(C\times C')/G$, which is effective.
Thus $tK_X$, and in turn $tK_S$ is effective. Since $K_S^2\geq K_X^2>0$,
it   follows from the Enriques classification of surfaces that
the minimal surface $S$ is of general type.
\qed

\medskip
To understand the remaining cases, and the geometry of the contraction $X\ra S$ as well,
we have to analyze the  two fibrations $\psi_i:X\ra\PP^1$, $i=1,2$ induced by the projections
$\pr_1:C\times C'\ra C$ and $\pr_2:C\times C'\ra C'$. Let us first record:

\begin{lemma}
\mylabel{geometrically connected}
We have $\psi_{i*}(\O_X)=\O_{\PP^1}$.
\end{lemma}

\proof
By the Main Theorem of Zariski, it suffices to check this at the generic point.
Let $E=k(C)$ be the function field of $C$, such that $L=E^G=k(\PP^1)$ 
is the function field of the projective line.
By construction, the generic fiber $X_\eta$ of $f_1:X\ra\PP^1$ is isomorphic to
$(C\otimes_k E)/G=(C_L\otimes_L E)/G$. In other words, the generic fiber 
is a twisted form of $C_L$ with respect to the \'etale topology. This description ensures
that $H^0(X_\eta,\O_{X_\eta})=L$, and the result follows.
\qed

\medskip
Now let $F_1\subset X$ and $F_2\subset X$ be the respective schematic fibers  for the   projections $\psi_1:X\ra\PP^1$
and $\psi_2:X\ra\PP^1$
containing the exceptional divisor $A_0+\sum_{i=0}^q\sum_{j=1}^{q-1}A_{i,j}$
for the resolution of singularities $X\ra(C\times C')/G$. 
Up to multiplicities, the fiber  $F_i$, $i=1,2$ is the union of this exceptional divisor $A_0+\sum_{i=0}^q\sum_{j=1}^{q-1}A_{i,j}$ 
with another integral curve $B_i\subset X$, which is birational to $\PP^1=C/G=C'/G$. 
One can say more about the curves $B_1,B_2\subset X$:

\begin{proposition}
\mylabel{minus-one curves}
The curves $B_1,B_2\subset X$ are $(-1)$-curves, and the
reduction of the fiber union $F_1\cup F_2\subset X$ has only simple  normal crossings.
\end{proposition}

\proof
By symmetry, it suffices to treat $B_1$. 
Being an integral component of a reducible fiber, it has selfintersection $B_1^2<0$.
Thus it suffices to check that $K_X\cdot B_1<0$. 
In light of Corollary \ref{canonical cycle} and Proposition \ref{fundamental cycle}, we have $K_{h}\cdot B_1\leq -(q-2)$.
The image of $qB_1$ on the normal surface $(C\times C')/G$ is a schematic fiber. 
Using the projection formula and computing intersections on $C\times C'$,
we deduce the value $h_*(B_1)\cdot K_{(C\times C')/G}= q-3$. The upshot is that 
$$
K_X\cdot B_1=K_{h}\cdot B_1 + K_{(C\times C')/G}\cdot h_*(B_1) \leq -1,
$$
whence $B_1$ is a $(-1)$-curve. Consequently, $K_h\cdot B_1=-(q-2)$.
In light of Corollary \ref{canonical cycle}, it follows that $B_i$ intersects the exceptional
divisor $A_0+\sum_{i=0}^q\sum_{j=1}^{q-1}A_{i,j}$ in precisely one component where the fundamental
cycle attains its minimal multiplicity, that is, $j=q-1$.
We conclude that the reduced fiber $F_{1,\red}$ is simple normal crossing.
\qed

\medskip
Let us  choose the indices $i$ for the $A_{i,j}$ so that $B_1\cdot A_{0,q-1}=B_2\cdot A_{q,q-1}=1$.
Thus the dual graph of the fiber union $F_1\cup F_2$ looks like this:

\vspace{2em}
\centerline{\includegraphics{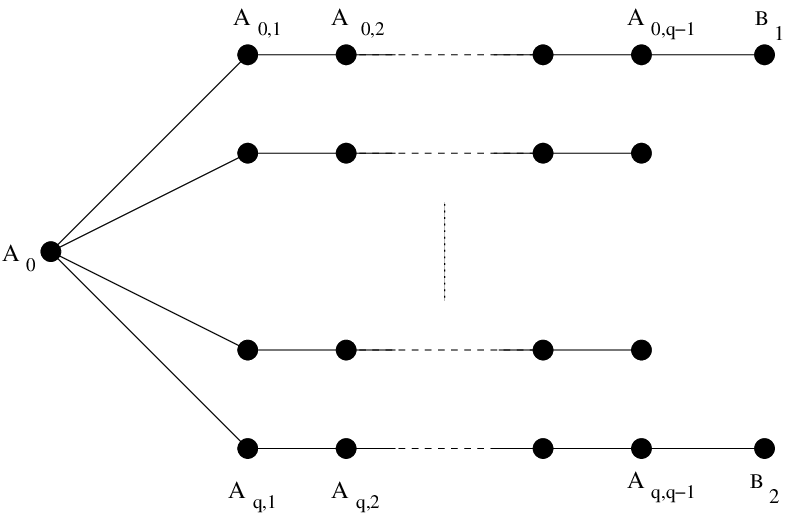}} 
\vspace{1em}
\centerline{Figure \stepcounter{figure}\arabic{figure}: Dual graph for fiber union $F_1\cup F_2$ on $X$.}
\vspace{1em}

We are now in position to compute the multiplicities occurring in the
schematic fibers $F_i\subset X$:

\begin{proposition}
\mylabel{schematic fibers}
The schematic fibers are given by  $F_1=Z+\sum_{j=1}^{q-1}jA_{0,j} + qB_1$ and $F_2=Z+\sum_{j=1}^{q-1}jA_{q,j} + qB_2$.
\end{proposition}

\proof
It suffices to treat $F_1$. By definition of the fundamental cycle $Z$, we have
$$
Z+\sum_{j=1}^{q-1}jA_{0,j} + qB_1= qA_0 +\sum_{i=1}^q\sum_{j=1}^{q-1} (q-j)A_{i,j} + \sum_{j=1}^{q-1}qA_{0,j} + q B_1,
$$
and it is a straightforward computation that this cycle is numerically trivial on all irreducible components
of $F_1$. It is therefore a rational multiple of $F_1$. Our cycle contains $B_1$ with multiplicity $q$.
But $qB_1$ is the strict transform
of the fiber for $\varphi_1:(C\times C')/G\ra\PP^1$,   by Proposition \ref{canonical quotient}.
Thus our cycle coincides with $F_1$.
\qed

\medskip
Note that the component $A_{0,q-1}$ has multiplicity one in the fiber $F_2\subset X$,
and similarly for $A_{q,q-1}\subset F_1\subset X$. In particular, the fibers are neither multiple nor wild.
Moreover:

\begin{corollary}
\mylabel{direct images}
The sheaf $R^1\psi_{i*}(\O_X)$ is locally free of rank $(q-1)(q-2)/2$,
and the formation of the direct image $\O_{\PP^1}=\psi_{i*}(\O_X)$ commutes with base change.
\end{corollary}

\proof
This follows from \cite{Raynaud 1970}, Theorem 7.2.1, because each geometric fiber of $\psi_i$ contains
a reduced irreducible component. 
\qed

\medskip
This leads to a very useful consequence concerning fundamental group:

\begin{corollary}
\mylabel{fundamental group}
The fundamental group $\pi_1(X)$ vanishes.
\end{corollary}

\proof
Let $X'\ra X$ be a finite \'etale covering with $X'$ nonempty.
We have to check that it has a section.
Consider the   projection  $\psi=\psi_1:X\ra\PP^1$, and let
$$
X'\lra T'\lra\PP^1
$$
be the Stein factorization for the composition $\psi':X'\ra\PP^1$, which is given by
$T'=\Spec(\psi'_*\O_{X'})$.
As explained in \cite{EGA IIIb}, Remark 7.8.10, the finite morphism $T'\ra\PP^1$ is \'etale,
since the equality $\O_{\PP^1}=\psi_*(\O_X)$ commutes with base change.
Moreover, the fiber $\psi^{-1}(0)\subset X$ is simply connected,  by the description in
Proposition \ref{schematic fibers}.
Using this as in the proof for 
\cite{SGA 1}, Expos\'e X, Theorem 1.3, we infer that the canonical map
$X'\ra X\times_{\PP^1}T'$, which is finite \'etale,  is actually an isomorphism.
But the projective line is simply connected, so
$T'\ra\PP^1$, whence also $X'\ra X$ has a section.
\qed

\medskip
For later use, we record:

\begin{proposition}
\mylabel{picard group}
The Picard group $\Pic(X)$ is a finitely generated free abelian group.
\end{proposition}

\proof
It already follows from Proposition \ref{irregularity} that $\Pic^0(X)=0$, whence $\Pic(X)=\NS(X)$ is finitely generated.
If $l\neq p$ is a prime different from the characteristic, then the elements $\shL\in\Pic(X)$ of order
$l$ yield nontrivial $\mu_l$-torsors $X'\ra X$, which are finite \'etale coverings of degree $l$.
Since the fundamental group vanishes, there are no such elements.

Finally, suppose there is an element $\shL$ of order $p$. To handle this case, first note
that the Artin--Schreier curve $C$ has $p$-rank $\sigma=0$. This follows from \cite{Subrao 1975}, Theorem 3.5,
applied to  $C\ra\PP^1$. In other words $H^1(C,\ZZ/p\ZZ)=0$, or equivalently, $\Pic(C)$ contains no element
of order $p$. Then the generic fiber $X_\eta$ for the second projection $\psi_2:X\ra\PP^1$, which is a twisted
form of $C$ over the function field of $\PP^1$, admits no invertible sheaf of order $p$.
Thus $\shL|X_\eta$ is trivial. It follows that $\shL\simeq\O_X(D)$, where $D$ is a divisor supported by
fibers of $\psi_2:X\ra\PP^1$. Now recall that the intersection form on the integral components of a given fiber
is negative semidefinite. Moreover, the radical is generated by the fiber, since the fibers have multiplicity one
by Proposition \ref{schematic fibers}.
Using that $\shL$ is numerically trivial, we first deduce that each $D_a$, $a\in\PP^1$ is a multiple of the   fiber,
and then that $D$ is linearly equivalent to the zero divisor, contradiction.
\qed

\medskip
The geometry of our surface $X$ simplifies further if we contract in a different way.
Indeed, the curves $B_1+\sum_{i=1}^{q-1}A_{0,j}$ and $B_2+\sum_{i=1}^{q-1}A_{q,j}$ are two disjoint 
\emph{exceptional curves of the first kind}, comprising altogether $2q$ irreducible components. 
Let 
$$
X\lra\tilde{S}
$$ 
be their contraction, such that $\tilde{S}$ is a smooth surface with
\begin{equation}
\label{change}
K_{\tilde{S}}^2=K_X^2 +2q.
\end{equation}
By choosing another minimal model $S$ if necessary, we may tacitly assume that the contraction $X\ra S$ factors over $\tilde{S}$.
Let $\tilde{A}_0,\tilde{A}_{i,j}\subset \tilde{S}$ be the images of the curves $A_0,A_{i,j}\subset X$ for $i\neq 0,q$. Then
$$
\tilde{Z}=\sum_{i=1}^{q-1}\sum_{j=1}^{q-1} (q-j)\tilde{A}_{i,j}
$$
is the image of the fundamental cycle $Z\subset X$ for the resolution of singularities $X\ra (C\times C')/G$.

\begin{proposition}
\mylabel{canonical divisor}
We have $K_{\tilde{S}}= (q-4)\tilde{Z}$, and $K_{\tilde{S}}^2=q(q-4)^2$.
\end{proposition}

\proof
Recall that $K_X^2= q(q^2-8q+14)$ by (\ref{selfintersection resolution}).
The map $X\ra\tilde{S}$ contracts successively $2q$ $(-1)$-curves,
such that
$$
K_{\tilde{S}}^2 =K_X^2 + 2q = q(q^2-8q+16)=q(q-4)^2.
$$
The canonical class $K_{\tilde{S}}$ is the image of the canonical class $K_X$.
Recall that we have 
$$
K_h=-(q-2)Z\quadand h^*K_{(C\times C')/G}= (q-3)(F_1+F_2).
$$
The latter coincides with the cycle
$2(q-3)Z$, up to components that are contracted by $X\ra \tilde{S}$. Consequently $K_{\tilde{S}}= (2(q-3) -(q-2))\tilde{Z}= (q-4)\tilde{Z}$.
\qed

\medskip
We now have an explicit description of the minimal model:

\begin{theorem}
\mylabel{minimal surface}
For $q\geq 4$, the surface $\tilde{S}$ is minimal, such that  $S=\tilde{S}$.
\end{theorem}

\proof
We have $\tilde{A}_0^2=2-q\neq -1$. Consequently, there is no $(-1)$-curve supported 
by $\tilde{Z}\subset\tilde{S}$. Since $K_{\tilde{S}}$ is effective, there is   no 
other $(-1)$-curve on $\tilde{S}$, and the result follows.
\qed

\medskip
From this we  easily determine the place  in the Enriques classification of surfaces.
Recall that a  \emph{weak del Pezzo surface} is a surface whose anticanonical divisor is nef and big.

\begin{corollary}
\mylabel{enriques classification}
The mimimal surface $S$ is of general type if $q\geq 5$,   a K3-surface 
for $q=4$, and   a  weak del Pezzo surface for $q=2,3$.
\end{corollary}

\proof
First suppose $q\geq 5$. Then $S=\tilde{S}$ is minimal, and $K_S$ is effective with $K_S^2>0$.
By the Enriques classification, $S$ is of general type.
For $q=4$, we  have $K_S=0$. Whence $S$ is either abelian, bielliptic, quasibielliptic, K3 or Enriques.
But $H^1(S,\O_S)=0$ by Proposition \ref{irregularity}, whence $S$ is either K3  or a classical Enriques
surface. In the latter case, $\pi_1(S)$ is cyclic of order two.
In light of Proposition \ref{fundamental group}, our $S$ must be a K3 surface.
Finally, suppose $q\leq 3$. Then $-K_{\tilde{S}}=(4-q)\tilde{Z}$ is effective,
and one easily checks that $\tilde{Z}\subset \tilde{S}$ is not an exceptional curve
of the first kind. Thus $-K_{\tilde{S}}$ is nef. It is also big, according to Proposition \ref{canonical divisor}.
The same necessarily holds for $-K_S$, hence $S$ is a weak del Pezzo surface.
\qed

\section{Canonical models and  canonical maps}
\label{Canonical model and canonical map}

In this   section we introduce  projective models for our surfaces.
If $q\geq 5$ then $S=\tilde{S}$ is a minimal surface of general type, and the homogeneous spectrum
$$\bar{S}=\Proj H^0(S,\bigoplus_{t\geq 0}\omega_S^{\otimes t})$$ is called
the \emph{canonical model} of $S$. The canonical morphism $S\ra\bar{S}$ is  
the contraction of all $(-2)$-curves, and the singularities on the normal surface $\bar{S}$ 
are at most rational double points. Clearly,
the integral curves $\tilde{A}_{i,j}\subset\tilde{S}=S$ are $(-2)$-curves, which get contracted.
It turns out that there are no more:

\begin{proposition}
\mylabel{canonical model}
For $q\geq5$, the canonical model $\bar{S}$ is obtained by contracting
the $(-2)$-curves $\tilde{A}_{i,j}$, $1\leq i,j\leq q-1$,
such that the singular locus of $\bar{S}$ comprises exactly $q-1$
rational double points of type $A_{q-1}$.
\end{proposition}

\proof
The exceptional curve for the contraction $\tilde{S}\ra\bar{S}$ is the union of all $(-2)$-curves, so
all   $\tilde{A}_{i,j}$, $1\leq i,j\leq q-1$ get contracted. Suppose there would be another
$(-2)$-curve   $\tilde{E}\subset\tilde{S}$, such that $\tilde{E}$ is disjoint from the canonical divisor $K_{\tilde{S}}=(q-4)\tilde{Z}$.
This implies that its strict transform $E\subset (C\times C')/G$ is contained in a smooth fiber
of one of the projections $\psi_i:(C\times C')/G\ra \PP^1$, contradiction.
\qed

\medskip
Note that the contraction $\tilde{S}\ra\bar{S}$ of the $(-2)$-curves $\tilde{A}_{i,j}$ makes sense for all prime powers $q=p^s$,
and it turns out that the normal surface $\bar{S}$ has a very satisfactory projective description.
To see this, let $\bar{Z}\subset\bar{S}$ be the Weil divisor defined as the image of $Z\subset X$.
Recall that the latter is nothing but the fundamental cycle for the resolution of singularities $X\ra(C\times C)/G$.
Then clearly
$$
\bar{Z}^2=q  \quadand \bar{Z}_\red=\PP^1\quadand \bar{Z}=q\bar{Z}_\red\quadand K_{\bar{S}}=(q-4)\bar{Z}.
$$
Since the local Picard group of a rational double point of type $A_{q-1}$ is cyclic of order $q$,
the Weil divisor $\bar{Z}\subset\bar{S}$ is actually Cartier.
The invertible sheaf $\bar{\shL}=\O_{\bar{S}}(\bar{Z})$ defines a rational map $\Phi_{\bar{Z}}:S\dashrightarrow\PP^n$,
with $n+1=h^0(\bar{\shL})$.

\begin{theorem}
The  invertible sheaf $\bar{\shL}=\O_{\bar{S}}(\bar{Z})$ is very ample, has $h^0(\bar{\shL})=4$,
and the image of the closed embedding $\Phi_{\bar{Z}}:\bar{S}\ra\PP^3$ is a normal surface
of degree $q$ whose singular locus consists of   $q-1$ rational double points
of type $A_{q-1}$, all lying on a line in $\PP^3$.
\end{theorem}

\proof
We first check that $\bar{\shL}$ is globally generated.
Let $\bar{F}_i\subset\bar{S}$, $i=1,2$ be the image of the
schematic fibers $F_i=\psi_i^{-1}(0)$ for the two projections
$\psi_i:X\ra\PP^1$. Using the multiplicities computed in Proposition \ref{schematic fibers}, we deduce $\bar{F}_1=\bar{Z}=\bar{F}_2$.
Now let $\bar{F}'_i\subset\bar{S}$ be the image of the schematic fibers
$F_i'=\psi_i^{-1}(\infty)$. Then $\bar{F}_i,\bar{F}'_i$ are linearly equivalent.
Using that $F'_1$ intersects $F_2$ only in  the  component $B_2$, we deduce that 
and $\bar{Z}\cap\bar{F}'_1$ is the image of the contracted curve $B_2$ under the canonical map $X\ra\bar{S}$.
By symmetry,  $\bar{Z}\cap\bar{F}'_2$ is the image of $B_1$. The upshot is that $\bar{Z}\cap\bar{F}'_1\cap\bar{F}'_2$ is empty,
hence $\bar{\shL}$ is globally generated. We also showed that $h^0(\bar{\shL})\geq 3$.

By the same ideas we verify that the resulting morphism $\Phi_{\bar{Z}}:\bar{S}\ra\PP^n$
is generically injective, where $n+1=h^0(\bar{\shL})$.
Let $x,x'\in X$ be two points not lying on the cycle $Z\cup F_1\cup F_2$ with $\psi_1(x)\neq\psi_1(x')$.
The   image on $\bar{S}$ of the schematic fiber $\psi_1^{-1}(\psi_1(x))\subset X$ is
a divisor linearly equivalent to $\bar{Z}$, which contains the image of $x$, but not  
the image of $x'$. Consequently, our map $\Phi_{\bar{Z}}$ is generically injective.

Next, we verify $h^0(\bar{\shL})=4$. To this end, consider the short exact sequence
$$
0\lra\O_{\bar{S}}\lra\bar{\shL}\lra\bar{\shL}|_{\bar{Z}}\lra 0.
$$
Since $h^1(\O_{\bar{S}})=0$, we have  $h^0(\bar{\shL})=h^0(\bar{\shL}|_{\bar{Z}})+1$.
The Weil divisors $i\bar{Z}_\red\subset\bar{S}$ yield short exact sequences
$$
0\lra\shK_i\lra \O_{(i+1)\bar{Z}_\red} \lra \O_{i\bar{Z}_\red}\lra 0.
$$
For $0\leq i<q$, the kernels are $\shK_i=\O_{\PP^1}(1-i)$, which follows from the intersection number
$\bar{Z}\cdot\bar{Z}_\red=1$ and
Proposition \ref{first blowing-up}, together with an application of \cite{Bayer; Eisenbud 1995}, Theorem 1.9. This gives the
estimate $h^0(\bar{\shL}_{\bar{Z}})\leq h^0\O_{\PP^1}(1)+h^0\O_{\PP^1}=3$.
If $h^0(\bar{\shL})=3$, then $\Phi_{\bar{Z}}:\bar{S}\ra\PP^2$
would be a finite surjective map of degree $q$, contradicting generic injectivity.
Thus $h^0(\bar{\shL})=4$.

Summing up, we have a morphism $\Phi_{\bar{Z}}:\bar{S}\ra\PP^3$ whose image $\hat{S}\subset\PP^3$
is a divisor of degree $q$, and the induced map $\nu:\bar{S}\ra\hat{S}$ is the normalization map.
Clearly, $\hat{S}$ is Cohen--Macaulay and Gorenstein. Using $\bar{\shL}=\Phi_{\bar{Z}}^*(\O_{\PP^3}(1))$
and $\omega_{\hat{S}}=\O_{\hat{S}}(q-4)$, we deduce $\omega_{\bar{S}}=\nu^*(\omega_{\hat{S}})$,
such that relative dualizing sheaf $\omega_{\bar{S}/\hat{S}}$ is trivial. From this it follows that the conductor locus
for the finite birational morphism $\nu:\bar{S}\ra\hat{S}$ is empty, such that $\nu$ is 
an isomorphism. Consequently $\bar{\shL}$ is very ample.
Finally observe that the 
image of $\bar{Z}_\red$, which contains the singular locus of $\bar{S}$, is a line, 
because $\bar{Z}\cdot\bar{Z}_\red=1$.
\qed

\medskip
Recall that a proper $k$-scheme $V_0$ is called \emph{liftable in the category of schemes}, if there exists 
a local ring $(R,\maxid_R)$ of characteristic zero with $k=R/\maxid_R$, together with a proper flat $R$-scheme $V$ with $V\otimes_W k=V_0$.
If such $V$ exists at least as an algebraic space, we say that $V_0$ is \emph{liftable in the category of algebraic spaces}.

\begin{corollary}
The surface $S$ is liftable in the category of algebraic spaces.
\end{corollary}

\proof
We first consider the projective normal surface $V_0=\bar{S}$.
Let $f_0\in k[X_0,\ldots,X_4]$ be a homogeneous polynomial of degree $\deg(f_0)=q$ so that 
$\bar{S}=V_+(f)$ as subscheme of $\PP^3=\Proj k[X_0,\ldots,X_4]$.
Choose   a homogeneous polynomial $f$ of degree $\deg(f)=q$ with coefficients in the ring of Witt vectors
$R=W(k)$ reducing to $f_0$ modulo $p$, and consider the proper $R$-scheme $V=V_+(f)\subset\PP^3_R$.
Then $V\ra\Spec(R)$ is flat and projective, such that $\bar{S}$ is liftable in the category of schemes.

According to a general result of Artin and Brieskorn \cite{Artin 1974},
there exists a finite extension $R\subset R'$ and a
simultaneous minimal resolution of singularities $W\ra V\otimes_RR'$.  Here, however, one knows only
that the total space $W$ is an algebraic space, although the individual fibers are projective. 
Thus $S=W_0$ is liftable in the category of algebraic spaces.
\qed

\section{Numerical invariants and geography}
\label{Numerical invariants and geography}

The \emph{Chern invariants} of a smooth proper surface are the numbers
$$
c_1^2=K^2\quadand c_2=e.
$$
They are paramount for minimal surfaces of general type, and
the study of occurrence and distribution of Chern invariants for minimal surfaces of general type is
referred to as \emph{surface geography}.

\begin{theorem}
\mylabel{chern invariants}
The Chern invariants for the smooth surface $\tilde{S}$ are given by the formulas
$$
K^2=q^3-8q^2+16q\quadand e=q^3-4q^2+6q,
$$
and   $\chi(\O_{\tilde{S}})=(q^3-6q^2+11q)/6 $.
\end{theorem}

\proof
We already computed  the value of $K^2$ in Proposition \ref{canonical divisor}.
Recall that $X\ra \tilde{S}$ is a sequence of $2q$ blowing-ups of $\tilde{S}$,
such that $e(X)=e(\tilde{S})+2q$.
To determine the Euler characteristic, we first examine the surface $X$
and its fibration $\psi_1:X\ra\PP^1$. According to Dolgachev's formula \cite{Dolgachev 1972}
we have  
$$
e(X) = e(X_{\bar{\eta}})e(\PP^1) + \sum (e(X_a)-e(X_{\bar{\eta}}) +\delta_a),
$$
where $X_{\bar{\eta}}$ is the geometric generic fiber,   the sum runs over all
closed point $a\in\PP^1$, and $\delta_a$ is \emph{Serre's measure of wild ramification}
attached to the Galois module $M_a=H^1(X_{\bar{\eta}},\ZZ/l\ZZ)$ at the point $a\in\PP^1$.
Here $l$ is any prime number different from $p$.
This   invariant is given by
$$
\delta_a = \sum_{i\geq 1} \frac{1}{[G:G_i]} \dim_{\FF_l}(M_a/M_a^{G_i}),
$$
Here $G$ is the Galois group of a finite Galois extension of the function field
$\kappa(\eta)$ trivializing the Galois module $M_a$, and $G_i\subset G$ are
the ramification subgroups for the induced extension of discrete valuation rings.
By the very construction, we may choose this extension induced from $C\times C'\ra(C\times C')/G$,
such that $G=\FF_q$. 

In light of (\ref{ramification groups}), we have $\delta_a=(q-1)\dim_{\FF_l}(M_a/M_a^G)$.
According to Proposition \ref{picard group}, the Picard group $\Pic(X)$ is finitely
generated.
Since the map $\Pic(X)\ra\Pic(X_\eta)$ is surjective, the group $\Pic(X_\eta)$
is finitely generated as well. Thus we may choose our prime $l$ so that $\Pic(X_\eta)$ contains
no nontrivial $l$-torsion. Since the Brauer group of the function field of $\PP^1$ vanishes
by Tsen's Theorem, it then follows that $M_a^G=0$, whence
$\dim_{\FF_l}(M_a/M_a^G)=(q-1)(q-2)$. In turn, we obtain the value for $e(\tilde{S})=e(X)-2q$.
Finally, Riemann--Roch for surfaces   $\chi=(K^2+e)/12$ yields the formula for $\chi(\O_{\tilde{S}})$.
\qed

\begin{remark}
There are two unfortunate misprints in \cite{Dolgachev 1972}: The correction
terms in Theorem 1.1 should appear with the sign lost in the proof while passing from equation (3.2) to (3.3).
In the definition of Serre's measure of wild ramification $\delta$ on top of page 305, the sum should   run only for $i\geq 1$,
which comes from the fact that the Swan character is the difference of the Artin character and
the augmentation character. Compare also \cite{Ogg 1967}.
\end{remark}

\medskip
We now can read off the \emph{geometric genus} $p_g=h^2(\O_X)=h^0(\omega_S)$:

\begin{corollary}
\mylabel{global geometric genus}
The  surface $S$ has geometric genus $p_g=(q^3-6q^2+11q-6)/6=(q-2)(q-3)(q-1)/6$.
\end{corollary}

\proof
According to Proposition \ref{irregularity},
we have $h^1(\O_S)=0$. The statement thus follows from $\chi(\O_S) =(q^3-6q^2+11q)/6$.
\qed

\medskip
Using our global invariants, we now can show that our bound on the geometric genus $p_g=\length R^1h_*(\O_X)$
of the singularity on $(C\times C')/G$ in Proposition \ref{geometric genus}
is actually an equality, at least in a special case:

\begin{corollary}
\mylabel{local geometric genus}
Suppose $q=p$ is prime. Then the singularity on $(C\times C')/G$ has geometric genus $p_g= p(p-1)(p-2)/6$.
\end{corollary}

\proof
The Leray--Serre spectral sequence for   $h:X\ra(C\times C')/G$ gives
an exact sequences
$$
H^1(X,\O_X)\lra H^0(R^1h_*(\O_X))\lra H^2(\O_{(C\times C')/G}) \lra H^2(\O_X)\lra 0.
$$
The term on the left vanishes by Proposition \ref{irregularity}.
We thus have
$$
p_g = h^0(\omega_{(C\times C')/G}) - h^2(\O_X).
$$
The first summand was computed in  Proposition \ref{invariant differentials},   the second
in Corollary \ref{global geometric genus}, and the result follows.
\qed

\begin{corollary}
\mylabel{plurigenera}
Suppose $q\geq 5$. Then the surface of general type $S$ has 
plurigenera 
$$
P_m = \frac{q^3-8q^2+16q}{2}(m^2-m) + \frac{q^3-6q^2+11q}{6}
$$ 
for all  $m\geq 2$.
\end{corollary}

\proof
We have $H^1(X,\omega_X^{\otimes m})= H^1(X,\omega^{\otimes(1-m)})=0$ for all $m\geq 2$ according to  
general results of Ekedahl  (\cite{Ekedahl 1988}, Theorem 1.7), and the value for 
$$
P_m=h^0(\omega_S^{\otimes m})=\chi(\omega_S^{\otimes m}) = \frac{m^2-m}{2} K_S^2 + \chi(\O_S),\quad m\geq 2
$$
follows from Riemann--Roch and the preceding  Theorem.
\qed

\begin{remark}
The \emph{Chern quotients}   of our surfaces $S$ asymptotically tend to
$$
\lim_{q\ra\infty} c_1^2/c_2 = 1.
$$
Note also that we always have $3c_2-c_1^2 = 2q(q-1)^2>0$,
such that the \emph{Bogomolov--Miyaoka--Yau inequality} $c_1^2<3c_2 $ holds.
Therefore our surfaces $S$ show
no exotic behavior with respect to geography.
\end{remark}

\begin{remark}
In the case $q=4$, we have $K_S=0$, and the preceding Theorem gives $e(S)=24$.
This yields another proof that $S$ is a K3 surface rather than an Enriques surface.
\end{remark}



\end{document}